\documentclass[11pt]{article}
\usepackage[left=1in,top=1in,right=1in,bottom=1in,head=.1in,nofoot]{geometry}

\usepackage{setspace,url,bm,amsmath} 

\usepackage{titlesec} 
\titlelabel{\thetitle.\quad} 

\usepackage{graphicx} 
\usepackage{bbm, amssymb}
\usepackage{latexsym}
\usepackage{caption}
\usepackage[margin=20pt]{subcaption}
\usepackage{hyperref}

\usepackage{color}

\usepackage[table]{xcolor}

\newcommand{\GG}[1]{}

\usepackage{amsthm}
\theoremstyle{definition}

\newtheorem*{theorem*}{Theorem}

\newtheorem*{corollary*}{Corollary}

\usepackage{natbib} 
\bibpunct{(}{)}{;}{a}{}{,} 

\usepackage{etoolbox} 
\apptocmd{\sloppy}{\hbadness 10000\relax}{}{} 

\usepackage{color}
\usepackage{listings}

\def\Var{\text{Var}}

\def\Cov{\text{Cov}}

\def\obs{\textnormal{obs}}

\def\tauS{\tau_\mathcal{S}}

\begin{document}
\doublespacing
\title{\bf Bridging Finite and Super Population Causal Inference}
\author{Peng Ding\footnote{Department of Statistics, University of California, Berkeley, California 94720 USA. Email: \url{pengdingpku@berkeley.edu}}, Xinran Li\footnote{Department of Statistics, Harvard University, Cambridge, Massachusetts 02138 USA. Email: 
\url{xinranli@fas.harvard.edu}}, and Luke W. Miratrix\footnote{Graduate School of Education and Department of Statistics, Harvard University, Cambridge, Massachusetts 02138 USA. Email: \url{luke_miratrix@gse.harvard.edu}}
}

\date{}
\maketitle

\begin{abstract}
There are two general views in causal analysis of experimental data: the super population view that the  units are an independent sample from some hypothetical infinite populations, and the finite population view that the potential outcomes of the experimental units are fixed and the randomness comes solely from the physical randomization of the treatment assignment. These two views differs conceptually and mathematically, resulting in different sampling variances of the usual difference-in-means estimator of the average causal effect. Practically, however, these two views result in identical variance estimators. By recalling a variance decomposition and exploiting a completeness-type argument, we establish a connection between these two views in completely randomized experiments. This alternative formulation could serve as a template for bridging finite and super population causal inference in other scenarios.

\end{abstract}

\section{Introduction}
\citet{Neyman:1923,neyman35} defined causal effects in terms of potential outcomes, and proposed an inferential framework viewing all potential outcomes of a finite population as fixed and the treatment assignment as the only source of randomness. This finite population view allows for easy interpretation free of any hypothetical data generating process of the outcomes, and is used in a variety of contexts \citep[e.g.,][]{Kempthorne:1952, hinkelmankempthorne, Copas:1973, rosenbaum2002observational, Imai:2008vp, freedmanOLSa, freedman2008randomization, Rosenbaum:2010book, aronow2013class, aronow2014sharp, abadie2014finite, miratrix2013adjusting, Ding:2014, lin.agnostic, middleton2015unbiased, ImbensRubinTextbook, chiba2015exact, rigdon2015randomization, li2016exact}. This approach is considered desirable because, in particular, it does not assume the data is somehow a representative sample of some larger (usually infinite) population.

Alternative approaches, also using the potential outcomes framework, assume that the potential outcomes are independent and identical draws from a hypothetical infinite population. Mathematical derivations under this approach are generally simpler, but the approach itself can be criticized because of this typically untenable sampling assumption. Furthermore, this approach appears to ignore the treatment assignment mechanism.

That being said, it is well known that the final variance formulae from either approach tend to be quite similar. For example, while the variance of the difference-in-means estimator for a treatment-control experiment under an infinite population model is different from the one under \citet{Neyman:1923}'s finite population formulation, this difference is easily represented as a function of the variance of the individual causal effects. Furthermore, this difference term is unidentifiable and is often assumed away under a constant causal effect model \citep{Neyman:1923, neyman35, hodges1970basic, rubin1990comment, reichardt1999justifying, FPP2007statistics}, or by appeals to the final estimators being ``conservative.''

Strictly speaking, the infinite population variance estimate gives a conservative (overly large) estimate of the finite population variance for the difference in means. As deriving infinite population variance expressions, relative to finite population variance expressions, tends to be more mathematically straightforward, we might naturally wonder if we could use infinite population expressions as conservative forms of finite population expressions more generally. In this work we show that in fact we can assume an infinite population model as an assumption of convenience, and derive formula from this perspective. This shows that we can thus consider the resulting formula as focused on the treatment assignment mechanism and not on a hypothetical sampling mechanism, i.e., we show variance derivations under the infinite population framework can be used as conservative estimators in a finite context.

Mathematically, this result comes from a variance decomposition and a completeness-style argument that characterizes the connection and the difference between these two views. The variance decomposition we use has previously appeared in \citet{Imai:2008vp}, \citet{ImbensRubinTextbook}, and \citet{balzer2016targeted}. 
The completeness-style argument, which we believe is novel in this domain, then sharpens the variance decomposition by moving from an expression on an overall average relationship to one that holds for any specific sample.
Our overall goal is simple: we wish to demonstrate that if one uses variance formula derived from assuming an infinite population sampling model, then the resulting inference one obtains will be correct with regards to the analogous samples-specific treatment effects (although it could be potentially conservative in that the standard errors may be overly large) regardless of the existence of any sampling mechanism. 

\section{Super Population, Finite Population, and Samples}

Assume that random variables $(Y(1), Y(0))$ represent the pair of potential outcomes of an infinite super population, from which we take an independent and identically distributed (IID) finite population of size $n$:
$$
\mathcal{S} =  \{  (Y_i(1), Y_i(0))  : i=1,\ldots, n\}.
$$
We first discuss completely randomized experiments, and comment on other experiments in Section \ref{sec::discussion}.
For a completely randomized experiment, we randomly assign $n_1$ units to receive treatment, leaving the remaining $n_0 = n-n_1$ units to receive control. Let $Z = (Z_1,\ldots,Z_n) \in \{0,1\}^n$ be the treatment assignment vector, which takes a particular value $(z_1,\ldots,z_n)\in \{0,1\}^n$ with probability $n_1!n_0!/n!$, for any $(z_1,\ldots,z_n)\in \{0,1\}^n$ with $\sum_{i=1}^n z_i=n_1$. 
The observed outcome for unit $i$ is then $Y_i^\obs = Z_i Y_i(1) + (1-Z_i) Y_i(0).$

At the super population level, the average potential outcomes are $E\{ Y(1) \}$ and $E\{  Y(0) \}$, and the average causal effect is
$$
  \tau = E\{ Y(1) - Y(0) \}.
$$
The population variances of the potential outcomes and individual causal effect are
$$
V_1 = \Var\{ Y(1) \},\quad
V_0 = \Var\{ Y(0)  \},\quad 
V_\tau =  \Var\{  Y(1) - Y(0) \}.
$$

At the finite population level, i.e., for a fixed sample $\mathcal{S}$, the average potential outcomes and average causal effect are 
$$
\bar{Y}_1 =  \frac{1}{n}  \sum_{i=1}^n   Y_i(1) ,\quad \bar{Y}_0 = \frac{1}{n} \sum_{i=1}^n Y_i(0) ,\quad
\tauS = \bar{Y}_1 - \bar{Y}_0  .
$$ 
The corresponding finite population variances of the potential outcomes and individual causal effects are
$$
S_1^2 = \frac{1}{n-1} \sum_{i=1}^n \{ Y_i(1) - \bar{Y}_1 \}^2,\quad
S_0^2 = \frac{1}{n-1} \sum_{i=1}^n \{ Y_i(0) - \bar{Y}_0 \}^2, \quad 
S_\tau^2 = \frac{1}{n-1} \sum_{i=1}^n (\tau_i - \tau)^2,
$$
where we use the divisor $n-1$ following the tradition of survey sampling \citep{Cochran:1977}.
All these quantities are fully dependent on $\mathcal{S}$.
In classical casual inference \citep{Neyman:1923}, the potential outcomes of these $n$ experimental units, $\mathcal{S}$, are treated as fixed numbers. Equivalently, we can consider such causal inference to be conducted conditional on $\mathcal{S}$ \citep[e.g.,][]{Copas:1973, rosenbaum2002observational, ImbensRubinTextbook, rigdon2015randomization, chiba2015exact}.

Regardless, we have two parameters---the population average treatment effect $\tau$, and the sample average treatment effect $\tauS$.
After collecting the data, we would want to draw inference about $\tau$ or $\tauS$.

Our primary statistics are the averages of the observed outcomes and the difference-in-means estimator:
$$
\bar{Y}_1^\obs =  \frac{1}{n_1 } \sum_{i=1}^n Z_i Y_i^\obs ,\quad \bar{Y}_0^\obs = \frac{1}{ n_0} \sum_{i=1}^n (1-Z_i) Y_i^\obs ,\quad
\widehat{\tau} =  \bar{Y}_1^\obs - \bar{Y}_0^\obs.
$$ 
We also observe the sample variances of the outcomes under treatment and control using 
$$
s_1^2 = \frac{1}{n_1-1} \sum_{i=1}^n Z_i  (   Y_i^\obs  - \bar{Y}_1^\obs ) ^2,\quad
s_0^2 = \frac{1}{n_0-1} \sum_{i=1}^n (1-Z_i) (   Y_i^\obs   - \bar{Y}_0^\obs ) ^2.
$$
We do not have the sample analogue of $S_\tau^2$ or $V_\tau$ because $Y_i(1)$ and $Y_i(0)$ are never jointly observed for any unit $i$ in the sample.

We summarize the infinite population, finite population and sample quantities in Table \ref{tb::three-levels}. 

\begin{table}[ht]
\centering
\caption{Means and variances at the super population, finite population and sample levels.}\label{tb::three-levels}
\begin{tabular}{|c|ccc|ccc|}
\hline 
&\multicolumn{3}{c|}{means} & \multicolumn{3}{c|}{variances} \\
\cline{2-4} \cline{5-7}
& treatment & control &  effect&treatment& control &  effect\\
\hline 
super population&$E\{ Y(1) \}$ & $E\{ Y(0) \}$ & $  \tau = E\{ Y(1) - Y(0)  \} $ & $V_1$ & $V_0$ & $V_\tau$ \\
finite population&$\bar{Y}_1$ & $\bar{Y}_0$ & $\tau_{\mathcal{S}} = \bar{Y}_1 - \bar{Y}_0 $ & $S_1^2$ & $S_0^2$ & $S_\tau^2$ \\
sample& $\bar{Y}_1^\obs$ & $\bar{Y}_0^\obs$ & $\widehat{\tau} = \bar{Y}_1^\obs- \bar{Y}_0^\obs$ & $s_1^2$ & $s_0^2$ & $-$\\
\hline 
\end{tabular}
\end{table}

\section{Deriving Complete Randomization Results with an Independent Sampling Model}
\label{sec::der}

The three levels of quantities in Table \ref{tb::three-levels} are connected via independent sampling and complete randomization. 
\citet{Neyman:1923}, without reference to any infinite population and by using the assignment mechanism as the only source of randomness, represented the assignment mechanism via an urn model, and found
\begin{eqnarray}
\Var(\widehat{\tau}\mid \mathcal{S}) = \frac{S_1^2}{n_1} + \frac{S_0^2}{n_0}  - \frac{S_\tau^2}{n}.
\label{eq::var-neyman}
\end{eqnarray} 
He then observed that the final term was unidentifiable but nonnegative, and thus if we dropped it we would obtain an upper bound of the estimator's uncertainty.

We next derive this result by assuming a hypothetical sampling mechanism from some assumed infinite super-population model of convenience.
This alternative derivation of the above result, which can be extended to other assignment mechanisms, shows how we can interpret formulae based on super-population derivations as conservative formulae for finite-sample inference.

\subsection{Sampling and randomization}
To begin, note that IID sampling of $\mathcal{S}$ from the super population implies three things: first, the finite population average causal effect $\tauS$ satisfies $E(\tauS) = \tau$, second
\begin{eqnarray}
\Var(\tauS) = 
 \Var\left[  \frac{1}{n} \sum_{i=1}^n\{  Y_i(1) - Y_i(0) \}  \right] 
=   \frac{V_\tau}{n},
\label{eq::iid}
\end{eqnarray}
and third, the sample variances are unbiased for the true variances:
\begin{eqnarray}
\label{eq::variances-unbiased}
E(S_1^2) = V_1,\quad
E(S_0^2) = V_0,\quad
E(S_\tau^2) = V_\tau.
\end{eqnarray}
Conditional on $\mathcal{S}$, randomization of the treatment $Z$ is the only source of randomness. 
In a completely randomized experiment, the outcomes in the treatment group form a simple random sample of size $n_1$ from $\{ Y_i(1): i=1,\ldots, n \}$, and the outcomes in the control group form a simple random sample of size $n_0$ from $\{  Y_i(0): i=1,\ldots, n \}$. 
Therefore, classical survey sampling theory \citep{Cochran:1977} for the sample mean and variance gives
\begin{eqnarray}
E(\bar{Y}_1^\obs\mid \mathcal{S}) = \bar{Y}_1,\quad
E(\bar{Y}_0^\obs\mid \mathcal{S}) = \bar{Y}_0,\quad 
E(\widehat{\tau} \mid \mathcal{S}) = \tau_{\mathcal{S}},\quad
E(s_1^2\mid \mathcal{S}) = S_1^2,\quad
E(s_0^2\mid \mathcal{S}) = S_0^2. 
\label{eq::cre}
\end{eqnarray}
We do not use the notation $E_{\mathcal{S}}(\cdot)$, because, depending on the contexts, such notation could either indicate expectation conditional on $\mathcal{S}$ or expectation averaged over $\mathcal{S}$. We therefore use conditional expectations and conditional variances explicitly.

If we do \emph{not} condition on $\mathcal{S}$, then the independence induced by the assignment mechanism means the outcomes under treatment are IID samples of $Y(1)$ and the outcomes under control are IID samples of $Y(0)$, and furthermore these samples are independent of each other.
This independence makes it straightforward to show that $\widehat{\tau}$ is unbiased for $\tau$ with super population variance
\begin{eqnarray}
\label{eq::var-super}
\Var(\widehat{\tau}) = \Var( \bar{Y}_1^\obs ) + \Var( \bar{Y}_0^\obs ) = \frac{V_1}{n_1} + \frac{V_0}{n_0}.
\end{eqnarray} 
This is the classic infinite population variance formula for the two sample difference-in-means statistic.
We could use it to obtain standard errors by plugging in $s_1^2$ and $s_0^2$ for the two variances.
This derivation, from \eqref{eq::iid}--\eqref{eq::var-super}, is motivated by the sampling assumption: the assignment mechanism makes the two samples independent, but it is the IID assumption and classic sampling theory that gives this result (along with asymptotical normality which allows for associated testing and confidence interval generation).

\subsection{Connecting the finite and infinite population inference with a variance decomposition}

We will now extend the above to indirectly derive the result on the variance of $\widehat{\tau}$ in finite population inference without explicitly enumerating the potential outcomes.
The variance decomposition formula implies
\begin{align}
\Var(\widehat{\tau}) &= E \{  \Var(\widehat{\tau}\mid \mathcal{S})  \}  + \Var\{ E (\widehat{\tau}\mid \mathcal{S}) \}  \label{eq::decompose}  \\
 &= E \{  \Var(\widehat{\tau}\mid \mathcal{S})  \}  + \Var(   \tauS   )   \nonumber \\
 &= E \{  \Var(\widehat{\tau}\mid \mathcal{S})  \}  +  \frac{V_\tau}{n} , \label{eq::proof}
\end{align} 
which further implies that the finite population variance of $\widehat{\tau}$ satisfies (using \eqref{eq::variances-unbiased})
\begin{eqnarray}
E \left\{  \Var(\widehat{\tau}\mid \mathcal{S})  \right\}  = \Var(\widehat{\tau}) -  \frac{V_\tau}{n}
= \frac{V_1}{n_1} + \frac{V_0}{n_0} -  \frac{V_\tau}{n}
=E \left\{  \frac{S_1^2}{n_1} + \frac{S_0^2}{n_0}  - \frac{S_\tau^2}{n}     \right\}.
\label{eq::variance}
\end{eqnarray}
Compare to the classic variance expression \eqref{eq::var-neyman}, which is this without the expectation.  
Here we have that \emph{on average} our classic variance expression holds.
Now, because this is true for any infinite population, as it is purely a consequence of the IID sampling mechanism and complete randomization, we can close the gap between \eqref{eq::var-neyman} and \eqref{eq::variance}. 
Informally speaking, because \eqref{eq::variance} holds as an average over many hypothetical super populations, it should also hold for any finite population at hand, and indeed it does, as we next show using a ``completeness'' concept from statistics \citep{lehmann2008testing}. 

\subsection{A ``Completeness'' argument}
First, define
\begin{align}
\label{eq::fS}
f(\mathcal{S}) \equiv \frac{S_1^2}{n_1} + \frac{S_0^2}{n_0} - \frac{S_\tau^2}{n} - \Var(\widehat{\tau}\mid \mathcal{S}),
\end{align} 
a function of a fixed finite sample $\mathcal{S}$, as the difference of the hypothesized finite sample variance formula and the actual finite sample variance.
Formula (\ref{eq::variance}) shows $E\{f(\mathcal{S})\}=0$. 
Now we are going to show the stronger result $f(\mathcal{S})=0$.

For any given sample $\mathcal{S}$, we have fixed sample quantities:
$
U =   \left(
Y_1(1)-\bar{Y}_1, \ldots, Y_n(1)-\bar{Y}_1\right),
$
and
$
W =  \left( Y_1(0)-\bar{Y}_0, \ldots, Y_n(0)-\bar{Y}_0 \right) .
$
Some algebra gives
\begin{align}\label{eq::var-tau}
\widehat{\tau}-\tau_{\mathcal{S}} = \frac{1}{n_1}\sum_{i=1}^n Z_iU_i-\frac{1}{n_0}\sum_{i=1}^n (1-Z_i)W_i.
\end{align}
According to \eqref{eq::var-tau}, the finite population variance $ \Var(\hat{\tau}\mid \mathcal{S}) =  E \{  ( \hat{\tau}-\tau_{\mathcal{S}})^2 \mid \mathcal{S} \} $ must be a quadratic form of $(U, W)$. Therefore, from \eqref{eq::fS},  $f(\mathcal{S})$ must be a quadratic form of $(U, W)$, and therefore must have the following expression:
\begin{align*}
f(\mathcal{S}) =  \sum_{i=1}^n \sum_{j=1}^n  a_{ij}U_iU_{j}+\sum_{i=1}^n \sum_{j=1}^n b_{ij}W_{i}W_j 
+ \sum_{i=1}^n  \sum_{j=1}^n c_{ij}U_{i}W_j.
\end{align*}
for some $a_{ij}$, $b_{ij}$ and $c_{ij}$.  
$f(\mathcal{S})$ is the same regardless of the ordering of units, so by symmetry the constants $a_{ij}$, $b_{ij}$ and $c_{ij}$ must be the same regardless of $i$ and $j$.  This gives
$
a_{ii} = a_{11},   a_{ij} = a_{12}, 
b_{ii} = b_{11},  b_{ij} = b_{12}, 
c_{ii} = c_{11},  c_{ij} = c_{12},
$
for all $i\neq j$.
Thus
\begin{align*}
f(\mathcal{S}) = &  a_{11}\sum_{i=1}^n U_j^2 + a_{12} \mathop{\sum\sum}_{i\neq j}   U_iU_j
+ b_{11}\sum_{i=1}^n W_i^2 + b_{12} \mathop{\sum\sum}_{i\neq j}  W_iW_j+
c_{11}\sum_{i=1}^n U_iW_i + c_{12} \mathop{\sum\sum}_{i\neq j} U_iW_j.
\end{align*}
Because $\sum_{i=1}^n U_i=0$, we have
\begin{align*}
 \sum_{i=1}^n U_i^2 + \mathop{\sum\sum}_{i\neq j} U_iU_j  = \left(  \sum_{i=1}^n  U_i \right)^2= 0.
\end{align*}
Similarly, as $\sum_{i=1}^n W_i = 0$, we have
\begin{align*}
\sum_{i=1}^n W_i^2 +\mathop{\sum\sum}_{i\neq j}  W_iW_j  = 0, \qquad 
\sum_{i=1}^n U_iW_i + \mathop{\sum\sum}_{i\neq j}  U_iW_j  =   \left(  \sum_{i=1}^n  U_i \right) \left(  \sum_{i=1}^n  W_i \right)  =  0.
\end{align*}
Use the above to replace our cross terms of $i \neq j$ with single index terms to write $f(\mathcal{S})$ as
\begin{align*}
f(\mathcal{S}) = &  a\sum_{i=1}^n U_i^2 
+ b\sum_{i=1}^n W_i^2 + 
c\sum_{i=1}^n U_iW_i,
\end{align*}
where $a=a_{11}-a_{12}, b=b_{11}-b_{12},$ and $c=c_{11}-c_{12}$. 
Because $E\{f(\mathcal{S})\}=0$, we have
$
aE(U_1^2)+ aE(W_1^2)+
cE(U_1W_1)
=0,
$
where
\begin{align*}
E(U_1^2) = \frac{n-1}{n}\Var\{Y(1)\}, \ 
E(W_1^2) = \frac{n-1}{n}\Var\{Y(0)\}, \ 
E(U_1W_1)^2 = \frac{n-1}{n}\Cov\{Y(1), Y(0)\}.
\end{align*}
Thus, 
\begin{align}\label{eq:constrain}
a\Var\{Y(1)\}+b\Var\{Y(0)\}+
c\Cov\{Y(1),Y(0)\}=0.
\end{align}
Because (\ref{eq:constrain}) holds for any populations regardless of its values of $\Var\{Y(1)\}, \Var\{Y(0)\}$ and $\Cov\{Y(1), Y(0)\}$, it must be true that 
$
a=b=c=0,
$
which implies $f(\mathcal{S})=0$.

\section{Discussion}\label{sec::discussion}

Equation~\eqref{eq::variance} relies on the assumption that the hypothetical infinite population exists, but Equation~\eqref{eq::var-neyman} does not. 
However, the completeness-style argument allowed us to make our sampling assumption only for convenience in order to prove \eqref{eq::var-neyman} by, in effect, dropping the expectation on both sides of \eqref{eq::variance}. 
Similar argument exists in the classical statistics literature; see \citet{efron1973stein} for the empirical Bayes view of Stein's estimator. 
While the final result is, of course, not new, we offer it as it gives an alternative derivation that does not rely on asymptotics such as a growing superpopulation or a focus on the properties of the treatment assignment mechanism.
 
Using  \citet{freedmanOLSa}'s results, \citet{aronow2014sharp} considered a super population with $N$ units, with the finite population being a simple random sample of size $n$. Letting $N\rightarrow \infty$, we can obtain similar results. 
We go in the other direction: we use the variance decomposition \eqref{eq::decompose} to derive the finite population variance from the super population one.    

This decomposition approach also holds for other types of experiments. 
First, for a stratified experiment, each stratum is essentially a completely randomized experiment. 
Apply the result to each stratum, and then average over all strata to obtain results for a stratified experiment. 
Second, because a matched-pair experiment is a special case of a stratified experiment with two units within each stratum, we can derive the Neyman-type variance \citep[cf.][]{Imai:2008vp, ImbensRubinTextbook} directly from that of a stratified experiment. 
Third, a cluster-randomized experiment is a completely randomized experiment on the clusters. If the causal parameters can be expressed as cluster-level outcomes, then the result can be straightforwardly applied \citep[cf.][]{aronow2013class,  middleton2015unbiased}. 
Fourth, for general experimental designs, the variance decomposition in \eqref{eq::decompose} still holds, and therefore we can modify the derivation of the finite population variance according to different forms of \eqref{eq::cre} and \eqref{eq::var-super}.

In a completely randomized experiment, the finite population sampling variance of $\widehat{\tau}$ in \eqref{eq::var-neyman} depends on three terms: the first two can be unbiasedly estimated by $s_1^2/n_1$ and $s_0^2/n_0$, but the third term $S_\tau^2/n$ is unidentifiable from the data. Assuming a constant causal effect model, $S_\tau^2=0$ and the variance estimators under both finite and super population inference coincide. However, \citet{Robins:1988}, \citet{aronow2014sharp} and \citet{Ding:2015jasa} demonstrated that the treatment variation term $S_\tau^2$ has a sharp lower bound that may be larger than $0$, which allows for more precise variance estimators under the finite population view. This demonstrates that we can indeed make better inference conditional on the sample we have. 
On the other hand, in this work, we showed that assuming an infinite population, while not necessarily giving the tightest  variance expressions, nonetheless gives valid (conservative) variance expressions from a finite-population perspective.
We offer this approach as a possible method of proof that could ease derivations for more complex designs.
More broadly, it is a step towards establishing that infinite population derivations for randomized experiments can be generally thought of as pertaining to their finite population analogs. 
Also see \citet{samii2012equivalencies} and \citet{lin.agnostic}, who provide alternative discussions of super population regression-based variance estimators under the finite population framework. 

Our discussion is based on the frequentists' repeated sampling evaluations of the difference-in-means estimator for the average causal effect. In contrast, \citet{Fisher:1935} proposed the randomization test against the sharp null hypothesis that $Y_i(1) = Y_i(0)$ for all units, which is numerically the same as the permutation test for exchangeable units sampled from an infinite population \citep{lehmann2008testing}. This connection becomes more apparent when only the ranks of the outcomes are used to construct the test statistics, as discussed extensively by \citet{Lehmann:1975}.

%
%
%

\section*{Acknowledgments}
We thank the Associate Editor, Dr. Peter Aronow, and three anonymous reviewers for helpful comments.

\bibliographystyle{apalike}
\bibliography{causalref}

\begin{thebibliography}{}

\bibitem[Abadie et~al., 2014]{abadie2014finite}
Abadie, A., Athey, S., Imbens, G.~W., and Wooldridge, J.~M. (2014).
\newblock Finite population causal standard errors.
\newblock Technical report, National Bureau of Economic Research.

\bibitem[Aronow et~al., 2014]{aronow2014sharp}
Aronow, P.~M., Green, D.~P., and Lee, D. K.~K. (2014).
\newblock Sharp bounds on the variance in randomized experiments.
\newblock {\em The Annals of Statistics}, 42:850--871.

\bibitem[Aronow and Middleton, 2013]{aronow2013class}
Aronow, P.~M. and Middleton, J.~A. (2013).
\newblock A class of unbiased estimators of the average treatment effect in
  randomized experiments.
\newblock {\em Journal of Causal Inference}, 1:135--154.

\bibitem[Balzer et~al., 2016]{balzer2016targeted}
Balzer, L.~B., Petersen, M.~L., and Laan, M.~J. (2016).
\newblock Targeted estimation and inference for the sample average treatment
  effect in trials with and without pair-matching.
\newblock {\em Statistics in Medicine}, 35:3717--3732.

\bibitem[Chiba, 2015]{chiba2015exact}
Chiba, Y. (2015).
\newblock Exact tests for the weak causal null hypothesis on a binary out come
  in randomized trials.
\newblock {\em Journal of Biometrics and Biostatistics},
  6:doi:10.4172/2155--6180.1000244.

\bibitem[Cochran, 1977]{Cochran:1977}
Cochran, W.~G. (1977).
\newblock {\em {Sampling Techniques}}.
\newblock John Wiley {\&} Sons, 3rd edition.

\bibitem[Copas, 1973]{Copas:1973}
Copas, J. (1973).
\newblock Randomization models for the matched and unmatched $2\times 2$
  tables.
\newblock {\em Biometrika}, 60:467--476.

\bibitem[Ding, 2016]{Ding:2014}
Ding, P. (2016).
\newblock A paradox from randomization-based causal inference (with
  discussion).
\newblock {\em Statistical Science}, in press.

\bibitem[Ding and Dasgupta, 2016]{Ding:2015jasa}
Ding, P. and Dasgupta, T. (2016).
\newblock A potential tale of two by two tables from completely randomized
  experiments.
\newblock {\em Journal of American Statistical Association}, 111:157--168.

\bibitem[Efron and Morris, 1973]{efron1973stein}
Efron, B. and Morris, C. (1973).
\newblock Stein's estimation rule and its competitors---an empirical bayes
  approach.
\newblock {\em Journal of the American Statistical Association}, 68:117--130.

\bibitem[Fisher, 1935]{Fisher:1935}
Fisher, R.~A. (1935).
\newblock {\em The {D}esign of {E}xperiments}.
\newblock Edinburgh, London: Oliver and Boyd, 1st edition.

\bibitem[Freedman et~al., 2007]{FPP2007statistics}
Freedman, D., Pisani, R., and Purves, R. (2007).
\newblock {\em Statistics}.
\newblock W. W. Norton \& Company, 4th edition.

\bibitem[Freedman, 2008a]{freedmanOLSa}
Freedman, D.~A. (2008a).
\newblock On regression adjustments in experiments with several treatments.
\newblock {\em The Annals of Applied Statistics}, 2:176--196.

\bibitem[Freedman, 2008b]{freedman2008randomization}
Freedman, D.~A. (2008b).
\newblock Randomization does not justify logistic regression.
\newblock {\em Statistical Science}, 23:237--249.

\bibitem[Hinkelmann and Kempthorne, 2008]{hinkelmankempthorne}
Hinkelmann, K. and Kempthorne, O. (2008).
\newblock {\em Design and Analysis of Experiments, Volume 1: Introduction to
  Experimental Design}.
\newblock New Jersey: John Wiley \& Sons, Inc., 2nd edition.

\bibitem[Hodges and Lehmann, 1964]{hodges1970basic}
Hodges, J. L.~J. and Lehmann, E.~L. (1964).
\newblock {\em Basic Concepts of Probability and Statistics}.
\newblock San Francisco: Holden-Day.

\bibitem[Imai, 2008]{Imai:2008vp}
Imai, K. (2008).
\newblock {Variance identification and efficiency analysis in randomized
  experiments under the matched-pair design}.
\newblock {\em Statistics in Medicine}, 27:4857--4873.

\bibitem[Imbens and Rubin, 2015]{ImbensRubinTextbook}
Imbens, G.~W. and Rubin, D.~B. (2015).
\newblock {\em Causal Inference for Statistics, Social and Biometrical
  Sciences: An Introduction}.
\newblock Cambridge: Cambridge University Press.

\bibitem[Kempthorne, 1952]{Kempthorne:1952}
Kempthorne, O. (1952).
\newblock {\em The {D}esign and {A}nalysis of {E}xperiments.}
\newblock New York: John Wiley and Sons.

\bibitem[Lehmann, 1975]{Lehmann:1975}
Lehmann, E.~L. (1975).
\newblock {\em Nonparametrics: Statistical Methods Based on Ranks}.
\newblock San Francisco: Holden-Day, 1st edition.

\bibitem[Lehmann and Romano, 2008]{lehmann2008testing}
Lehmann, E.~L. and Romano, J.~P. (2008).
\newblock {\em Testing Statistical Hypotheses}.
\newblock New York : Wiley, 3rd edition.

\bibitem[Li and Ding, 2016]{li2016exact}
Li, X. and Ding, P. (2016).
\newblock Exact confidence intervals for the average causal effect on a binary
  outcome.
\newblock {\em Statistics in Medicine}, 35:957--960.

\bibitem[Lin, 2013]{lin.agnostic}
Lin, W. (2013).
\newblock Agnostic notes on regression adjustments to experimental data:
  Reexamining {F}reedman's critique.
\newblock {\em The Annals of Applied Statistics}, 7:295--318.

\bibitem[Middleton and Aronow, 2015]{middleton2015unbiased}
Middleton, J.~A. and Aronow, P.~M. (2015).
\newblock Unbiased estimation of the average treatment effect in
  cluster-randomized experiments.
\newblock {\em Statistics, Politics and Policy}, 6:39--75.

\bibitem[Miratrix et~al., 2013]{miratrix2013adjusting}
Miratrix, L.~W., Sekhon, J.~S., and Yu, B. (2013).
\newblock {Adjusting treatment effect estimates by post-stratification in
  randomized experiments}.
\newblock {\em Journal of the Royal Statistical Society: Series B (Statistical
  Methodology)}, 75:369--396.

\bibitem[Neyman, 1923]{Neyman:1923}
Neyman, J. (1923).
\newblock On the application of probability theory to agricultural experiments.
  essay on principles (with discussion). section 9 (translated). reprinted ed.
\newblock {\em Statistical Science}, 5:465--472.

\bibitem[Neyman, 1935]{neyman35}
Neyman, J. (1935).
\newblock Statistical problems in agricultural experimentation (with
  discussion).
\newblock {\em Supplement to the Journal of the Royal Statistical Society},
  2:107--180.

\bibitem[Reichardt and Gollob, 1999]{reichardt1999justifying}
Reichardt, C.~S. and Gollob, H.~F. (1999).
\newblock Justifying the use and increasing the power of a $t$ test for a
  randomized experiment with a convenience sample.
\newblock {\em Psychological Methods}, 4:117--128.

\bibitem[Rigdon and Hudgens, 2015]{rigdon2015randomization}
Rigdon, J. and Hudgens, M.~G. (2015).
\newblock Randomization inference for treatment effects on a binary outcome.
\newblock {\em Statistics in Medicine}, 34:924--935.

\bibitem[Robins, 1988]{Robins:1988}
Robins, J.~M. (1988).
\newblock Confidence intervals for causal parameters.
\newblock {\em Statistics in Medicine}, 7:773--785.

\bibitem[Rosenbaum, 2002]{rosenbaum2002observational}
Rosenbaum, P.~R. (2002).
\newblock {\em Observational Studies}.
\newblock New York: Springer, 2nd edition.

\bibitem[Rosenbaum, 2010]{Rosenbaum:2010book}
Rosenbaum, P.~R. (2010).
\newblock {\em {Design of Observational Studies}}.
\newblock New York: Springer.

\bibitem[Rubin, 1990]{rubin1990comment}
Rubin, D.~B. (1990).
\newblock Comment: Neyman (1923) and causal inference in experiments and
  observational studies.
\newblock {\em Statistical Science}, 5:472--480.

\bibitem[Samii and Aronow, 2012]{samii2012equivalencies}
Samii, C. and Aronow, P.~M. (2012).
\newblock On equivalencies between design-based and regression-based variance
  estimators for randomized experiments.
\newblock {\em Statistics and Probability Letters}, 82:365--370.

\end{thebibliography}

\end{document}